\documentclass[12pt]{article}

\usepackage{amsfonts}
\usepackage{graphicx}

\newcommand{\R}{{\mathbb R}}
\newcommand{\N}{{\mathbb N}}
\newcommand{\Z}{{\mathbb Z}}
\newcommand{\C}{\mathbb{C}}
\newcommand{\F}{\mathbb{F}}

\newcommand{\ad}{{\rm ad}}
\newcommand{\Ad}{{\rm Ad}}

\newcommand{\Int}{{\rm Int}}
\newcommand{\Gl}{{\rm Gl}}
\newcommand{\gl}{{\frak{gl}}}
\newcommand{\tr}{{\rm tr}}

\renewcommand{\int}{{\rm int} \,}

\newcommand{\g}{\mathfrak{g}}

\newcommand{\n}{\mathfrak{n}}

\newcommand{\ov}[1]{{\overline{#1}}}

\newcommand{\wh}[1]{{\widehat{#1}}}
\newcommand{\wt}[1]{{\widetilde{#1}}}

\newtheorem{teorema}{Theorem}[section]
\newtheorem{lema}[teorema]{Lemma}
\newtheorem{corolario}[teorema]{Corollary}

\newenvironment{prova}{\noindent {\bf Proof:}}{\hfill $\qed $ \newline}

\newenvironment{exemplo}{\noindent {\bf Example:}}{}
\def\squarebox#1{\hbox to #1{\hfill\vbox to #1{\vfill}}}
\newcommand{\qed}{\hspace*{\fill}
\vbox{\hrule\hbox{\vrule\squarebox{.667em}\vrule}\hrule}\smallskip}

\newcommand{\paragrafo}{\vrule height 0pt width 0pt depth 0pt\hbox to\parindent{\hss}}

\begin{document}

\title{A note on the Jordan decomposition}
\author{Mauro Patr\~{a}o\footnote{Supported by FAPEDF grant
no.??/????},\, La{\'e}rcio Santos\footnote{Supported by FAPESP grant
no.??/????}\, and Lucas Seco\footnote{Supported by FAPESP grant
no.??/????}\\} \maketitle

\begin{abstract}
In this article we prove that the elliptic, hyperbolic and nilpotent
(or unipotent) additive (or multiplicative) Jordan components of an
endomorphism $X$ (or an isomorphism $g$) of a finite dimensional
vector space are given by polynomials in $X$ (or in $g$). By using
this, we provide a simple proof that, for an element $X$ of a linear
semisimple Lie algebra $\g$ (or $g$ of a linear semisimple connected
Lie group $G$), its three Jordan components lie again in the algebra
(in the group). This was previously unknown for linear Lie groups
other then $\Int(\g)$. This implies that, for this class of algebras
and groups, the usual linear Jordan decomposition coincides with the
abstract Jordan decomposition.
\end{abstract}

\section{Introduction}
Let $V$ be a finite dimensional real vector space and $T$ a linear
map of $V$. The most usual Jordan decomposition writes $T$ as a
commuting sum of a semisimple and a nilpotent maps.  They are
called the semisimple and nilpotent additive Jordan components of
$T$ and are given as polynomials in $T$ (see Theorem 13, p.267 of
\cite{hk} or Proposition 4.2, p.17 of \cite{hum}). One
can go further and write the semisimple component as a commuting
sum of an elliptic and a hyperbolic components which also commutes
with the nilpotent component. When $T$ is invertible, there is an
analogous multiplicative Jordan decomposition which writes $T$ as
a commuting product of an elliptic, a hyperbolic and a unipotent
components (see Section IX.7, p.430 of \cite{helgason}). In
Section \ref{sec1}, our main results show that the elliptic and
hyperbolic components of both additive and multiplicative Jordan
decomposition of $T$ are given as polynomials in $T$ and the same
happens for the unipotent component.

One can extend these decompositions to the context of semisimple
Lie algebras and groups in the following manner (see
\cite{helgason,var,warner}). Let $\g$ be a semisimple Lie algebra
and $G$ be a Lie group with Lie algebra $\g$. Let $\ad: \g \to
\gl(\g)$ be the adjoint representation of $\g$ and let $\Ad: G \to
\Gl(\g)$ be the adjoint representation of $G$. For $X \in \g$ we
say that $X = E + H + N$ is an abstract Jordan decomposition of
$X$ if $E, H, N \in \g$ commute, $\ad(E)$ is additively elliptic,
$\ad(H)$ is additively hyperbolic and $\ad(N)$ is nilpotent.  For
$g \in G$ we say that $g = ehu$ is an abstract Jordan
decomposition of $g$ if $e, h, u \in G$ commute, $\Ad(e)$ is
elliptic, $\Ad(h)$ is hyperbolic and $\Ad(u)$ is unipotent. In
Section \ref{sec2}, by using the results of Section \ref{sec1}, we
provide a simple proof that, for an element $X$ of a linear
semisimple Lie algebra $\g$ (or $g$ of a linear semisimple
connected Lie group $G$), its three Jordan components lie again in
the algebra (in the group). This was previously unknown for linear
Lie groups other then $\Int(\g)$. This implies that, for this
class of algebras and groups, the usual linear Jordan
decomposition coincides with the abstract Jordan decomposition.

We now introduce some preliminary definitions and notations.
Defining the complex vector space $V_\C = \{ u + iv:\, u,v \in
V\}$ we have that $V \subset V_\C$. For $X \in \gl (V)$ we put
$X(u + iv) = Xu + iXv$ so that $\gl(V) \subset \gl(V_\C)$. Since
the determinant of $g \in \Gl(V)$ seen as an operator of $V$ or
$V_\C$ coincide we also have that $\Gl(V) \subset \Gl(V_\C)$. Let
$X \in \gl(V)$. As usual, we say that $X$ is semisimple if it is
diagonalizable in $V_\C$ and that $X$ is nilpotent if there exists
$n \in \N$ such that $X^n = 0$. We say that $X$ is elliptic
(hyperbolic) in the additive case if it is semisimple and its
eigenvalues are purely imaginary (real). Now let $g \in \Gl(V)$.
We say that $g$ is elliptic (hyperbolic) in the multiplicative
case if it is semisimple and its eigenvalues have absolute value equal to one (are
real positive). We say that $g$ is unipotent if $g - I$ is
nilpotent. The proof of the following result is straightforward.

\begin{lema}\label{lemacomutacao}
Let $L, T$ be two commuting linear maps of $V$.
\begin{itemize}
  \item[1.] If both $L, T$ are semisimple, elliptic or hyperbolic, then
  $L + T$ in the additive case (or $LT$ in the multiplicative case) is semisimple, elliptic or hyperbolic.

  \item[2.] If both $L, T$ are nilpotent (or unipotent), then
  $L + T$ is nilpotent (or $LT$ is unipotent).

  \item[3.] If $T$ is simultaneously semisimple and nilpotent (or
  semisimple and unipotent) then $T=0$ (or $T=I$).

  \item[4.]  If $T$ is elliptic and hyperbolic, then $T=0$ in the additive case (or $T=I$ in the multiplicative case).
\end{itemize}
\end{lema}
%
%\begin{prova}
%Item (1) follows using the definitions, since $L$ and $T$ commute
%and thus they have the same eigenvectors. For item (2), first assume
%that $L^n = T^n = 0$, for some $n \in \N$. We have that $(L+T)^{2n}$
%vanishes, since its binomial coefficients are given by ${2n \choose
%k} L^kT^{2n -k } = 0$. Now assume that $L = I + N_L$ and $T = I +
%N_T$, where $N_L, N_T$ are nilpotent. Since $L, T$ commute, it
%follows that $N_L$, $N_T$ also commute. Therefore $LT = I + N_L +
%N_T + N_LN_T$ is unipotent, since $N_L + N_T + N_LN_T$ is a
%commutative sum of nilpotent maps. For item (3),  $T$
%
%\end{prova}

\section{General linear algebra and group}\label{sec1}

Denote by $\F[x]$ the ring of the polynomials in $x$ with
coefficients in $\F = \C \mbox{ or } \R$. We denote by $\ov{p}(x)$
the polynomial whose coefficients are the conjugate of the
coefficients of $p(x) \in \F[x]$. Thus $p(x) \in \R[x]$ if and only
if $\ov{p}(x) = p(x)$. Let $T \in \gl(V)$ and consider the following
ring homomorphism
\[
\F[x] \to \gl(V), \quad p(x) \mapsto p(T),
\]
where $p(x) = a_0 + a_1x + \cdots + a_m x^m$ and $p(T) =  a_0I +
a_1T + \cdots + a_m T^m$. We denote by $\F(T)$ the image of this
homomorphism. From now on, we will denote simply by $p$ both $p(x)$
and $p(T)$. It will be clear from the context which polynomial is
considered. The kernel of the above homomorphism is the principal
ideal generated by $p_T$, the so called minimal polynomial of $T$.
Since $T \in \gl(V)$, it follows that $p_T \in \R[x]$.  Since $p_T
\in \R[x]$ we can factor it over $\C$ as
\[
p_T = p_1 \ov{p}_1 \cdots p_l \ov{p}_l p_{l+1} \cdots p_n
\]
where $p_k(x) = (x - \lambda_k)^{m_k}$, $\lambda_k$ has imaginary
part for $k=1,\ldots,l$ and $\lambda_k$ is real for $k =
l+1,\ldots,n$.  Note that, since the characteristic polynomial of
$T$ divides $p_T(x)$, we have that the eigenvalues of $T$ are
$\lambda_k, \, \ov{\lambda_k}$, for $k = 1,\ldots,n$.

\begin{lema}\label{lemaprojecoes}
There exist polynomials $\pi_k$, where $1 \leq k \leq n$, and $l$
such that
\begin{itemize}
\item[(i)] $\pi_k \in \C[x]$, for $1 \leq k \leq l$, and $\pi_k \in
\R[x]$, for $l+1 \leq k \leq n$.

\item[(ii)] If $r \neq s$ then $\pi_r \pi_s$, $\ov{\pi}_r \pi_s$ and
$(x - \lambda_r)^{m_r}\pi_r$ are multiples of $p_T$. We also have that
$\ov{\pi}_r \pi_r$ is multiple of $p_T$ for $r=1,\ldots,l$.

\item[(iii)] $1 = \sum_{k=1}^l( \pi_k + \ov{\pi}_k ) + \sum_{k=l+1}^n
\pi_k$.
\end{itemize}
\end{lema}
\begin{prova}
For $1 \leq k \leq l$, we define the polynomials
\[
q_k = p_1 \ov{p}_1 \cdots p_k \wh{\ov{p}}_k \cdots p_l \ov{p}_l
p_{l+1} \cdots p_n
\]
whose conjugates are given by
\[
\ov{q}_k = p_1 \ov{p}_1 \cdots \wh{p}_k \ov{p}_k \cdots p_l \ov{p}_l
p_{l+1} \cdots p_n,
\]
where the factor below $\,\, \wh{} \,\,$ is omitted. For $l+1 \leq k
\leq n$, we define the polynomials
\[
q_k = p_1 \ov{p}_1 \cdots p_l \ov{p}_l p_{l+1} \cdots \wh{p}_k
\cdots p_n.
\]
Since the constant polynomials are the only polynomials dividing all
$q_k$ and $\ov{q}_k$, where $1 \leq k \leq n$, it follows that the
ideal generated by them is all of $\C[x]$. Thus there exists
polynomials $a_k$, $b_k$ and $c_k$ such that
\[
1 = \sum_{k=1}^l (a_kq_k + b_k\ov{q}_k) + \sum_{k=l+1}^n c_kq_k.
\]
Adding the above equation with its conjugate and dividing by two, we
can assume that $b_k = \ov{a}_k$ and $c_k = \ov{c}_k$.  Defining the
polynomials $\pi_k = a_k q_k$, $k = 1,\ldots,l$, and $\pi_k = b_k
q_k$, $k = l+1,\ldots,n$ we obtain the result.
\end{prova}

\begin{lema}\label{lemapolT}
Applying the above polynomials to $T$ we have the following.
\begin{itemize}
\item[(i)] $I = \sum_{k=1}^l( \pi_k + \ov{\pi}_k ) +
\sum_{k=l+1}^n \pi_k$.

\item[(ii)] If $r \neq s$, then $\pi_r \pi_s = 0$, $\ov{\pi}_r
\pi_s = 0$ and $(T - \lambda_r)^{m_r}\pi_r = 0$. We also have that
$\ov{\pi}_r \pi_r = 0$, for $r=1,\ldots,l$.
%$\ov{\pi}_r \pi_s = 0$, $(T - \lambda_r)^{m_r}\pi_r = 0$ and if $r \neq s$ then $\pi_r \pi_s = 0$.

\item[(iii)] For $r=1,\ldots,n$ we have $\pi_r^2 = \pi_r$.
\end{itemize}
\end{lema}
\begin{prova}
Items (i) and (ii) are immediate from the previous lemma and the
definition of the minimal polynomial.  For item (iii), apply $\pi_k$
to both sides of item (i) and use item (ii).
% The last item is
% immediate from items (ii) and (iii).
\end{prova}

Now we make the following remarks. By items (i) and (ii) of the
above lemma we have that $V_\C$ is the direct sum of the images of
the projections $\pi_k$, $\ov{\pi}_k$. Also, let
$$
A_r = \sum_{k=1}^l (a_{rk} \pi_k + b_{rk} \ov{\pi}_k) +
\sum_{k=l+1}^n c_{rk} \pi_k,
$$
where $r =1,2$ and $a_{rk}$, $b_{rk}$, $c_{rk}$ are polynomials in
$T$.  Again by the above lemma, we have that
\[
A_1 A_2 = \sum_{k=1}^l (a_{1k} a_{2k} \pi_k + b_{1k} b_{2k}
\ov{\pi}_k) + \sum_{k=l+1}^n c_{1k} c_{2k} \pi_k.
\]
We now obtain the description of the additive Jordan components as
polynomials.

\begin{teorema}\label{teojordanaditivo}
Let $X \in \gl(V)$ and put $\lambda_k = u_k + i v_k$. Then $X$ can
be written uniquely as a commutative sum $X = E + H + N$, where
$E$ is elliptic, $H$ is hyperbolic and $N$ is nilpotent and they
are given by the following real polynomials
\begin{itemize}
\item[(i)] $E = \sum_{k=1}^l( iv_k\pi_k + \ov{iv_k \pi}_k )$.

\item[(ii)] $H = \sum_{k=1}^l u_k(\pi_k + \ov{\pi}_k ) +
\sum_{k=l+1}^n u_k\pi_k$.

\item[(iii)] $N = \sum_{k=1}^l( (X - \lambda_k)\pi_k + (X - \ov{\lambda}_k)\ov{\pi}_k ) +
\sum_{k=l+1}^n (X - \lambda_k)\pi_k$.
\end{itemize}
\end{teorema}
\begin{prova}
Let $S = \sum_{k=1}^l( \lambda_k \pi_k + \ov{\lambda_k \pi}_k ) +
\sum_{k=l+1}^n \lambda_k \pi_k$.  Note that $N = X - S$ and that $S
= E + H$ and thus $X = E + H + N$.  Using the remarks after Lemma
\ref{lemapolT}, it is immediate that $S$ is semisimple, $E$ is
elliptic, $H$ is hyperbolic and that
$$
N^m = \sum_{k=1}^l( (X - \lambda_k)^m\pi_k + (X -
\ov{\lambda}_k)^m\ov{\pi}_k ) + \sum_{k=l+1}^n (X -
\lambda_k)^m\pi_k.
$$
Taking $m = \max_k\{ m_k \}$, by item (ii) of Lemma \ref{lemapolT},
it follows that $N^m = 0$.

For the uniqueness, consider the commuting sum $X = \wt{E} + \wt{H}
+ \wt{N}$, with $\wt{E}$ elliptic, $\wt{H}$ hyperbolic and $\wt{N}$
nilpotent. Define $\wt{S} = \wt{E} + \wt{H}$.  Since $\wt{E}$ and
$\wt{H}$ commute, by Lemma \ref{lemacomutacao}, we have that
$\wt{S}$ is semisimple. Since $E$, $H$, $N$ are polynomials in $X$,
they commute with $\wt{E}$, $\wt{H}$, $\wt{N}$.  Using that $X = S +
N = \wt{S} + \wt{N}$, by Lemma \ref{lemacomutacao}, we have that $S
- \wt{S} = \wt{N} - N$ is both semisimple and nilpotent and thus $S
= \wt{S}$ and $\wt{N} = N$.  Now using that $S = E + H = \wt{E} +
\wt{H}$,  by Lemma \ref{lemacomutacao}, we have that $E - \wt{E} =
\wt{H} - H$ is both elliptic and hyperbolic and thus $E = \wt{E}$
and $\wt{H} = H$.
\end{prova}

The following result provides the description of the
multiplicative Jordan components as polynomials.

\begin{teorema}\label{teojordanmult}
Let $g \in \Gl(V)$ and put $\lambda_k = u_k + i v_k$. Then $g$ can
be written uniquely as a commutative product $g = ehu$, where $e$
is elliptic, $h$ is hyperbolic and $u$ is unipotent and they are
given by the following real polynomials
\begin{itemize}
\item[(i)] $e = \sum_{k=1}^l |\lambda_k|^{-1}(\lambda_k\pi_k + \ov{\lambda_k \pi}_k ) +
\sum_{k=l+1}^n |\lambda_k|^{-1} \lambda_k \pi_k$.

\item[(ii)] $h = \sum_{k=1}^l |\lambda_k|(\pi_k + \ov{\pi}_k ) +
\sum_{k=l+1}^n |\lambda_k|\pi_k$.

\item[(iii)] $u = I + N\left(\sum_{k=1}^l( \lambda_k^{-1}\pi_k + \ov{\lambda_k^{-1}\pi}_k ) +
\sum_{k=l+1}^n \lambda_k^{-1}\pi_k\right)$,
\end{itemize}
where $N$ is the nilpotent component of $g$.  Furthermore, we have
that $h = {\rm e}^{H}$, where
\[
H = \sum_{k=1}^l \log(|\lambda_k|)(\pi_k + \ov{\pi}_k ) +
\sum_{k=l+1}^n \log(|\lambda_k|)\pi_k.
\]
\end{teorema}
\begin{prova}
By the proof of Theorem \ref{teojordanaditivo} we have that $g = S +
N$.  Noting that $u = I + NS^{-1}$, we have that $Su = S + N = g$.
It is immediate that $S = eh$, and thus $g = ehu$.  Since $S$, $N$
commute and $N$ is nilpotent, it follows that $u$ is unipotent.
Using the remarks after Lemma \ref{lemapolT}, it is immediate that
$e$ is elliptic, that $h$ is hyperbolic and that $h = {\rm e}^{H}$.

For the uniqueness, consider the commuting product $g = \wt{e}
\wt{h} \wt{u}$, with $\wt{e}$ elliptic, $\wt{h}$ hyperbolic and
$\wt{u}$ unipotent. Define $\wt{S} = \wt{e}   \wt{h}$.  Since
$\wt{e}$ and $\wt{h}$ commute, by Lemma \ref{lemacomutacao}, we have
that $\wt{S}$ is semisimple. Since $e$, $h$, $u$ are polynomials in
$g$, they commute with $\wt{e}$, $\wt{h}$, $\wt{u}$.  Using that $g
= Su    = \wt{S}   \wt{u}$, by Lemma \ref{lemacomutacao}, we have
that $\wt{S}^{-1}S= \wt{u}u^{-1}$ is both semisimple and unipotent
and thus $S = \wt{S}$ and $\wt{u} = u$.  Now using that $S = eh =
\wt{e}\wt{h}$, by Lemma \ref{lemacomutacao}, we have that
$\wt{e}\,^{-1}e = \wt{h}h^{-1}$ is both elliptic and hyperbolic and
thus $e = \wt{e}$ and $\wt{h} = h$.
\end{prova}

\begin{exemplo}
In this example we display polynomials which give the additive and
multiplicative Jordan components of the linear map
$$
T = \left(
\begin {array}{cccc}
1&1&0&0\\\noalign{\medskip}-1&1&0&0\\\noalign{\medskip}0&0&2&1\\\noalign{\medskip}0&0&0&2
\end{array}
\right),
$$
which has $p_T(x) = (x - (1 + i))(x - (1 - i))(x-2)^2$. By using
the Euclidean algorithm, we obtain the following polynomials which
satisfy Lemma \ref{lemaprojecoes}
$$
\pi_1(x) = \frac{1}{4}\, \left( x-1+i \right)  \left( x-2 \right)
^{2},\quad
\pi_2(x) = -\frac{1}{2}\, \left( x-3 \right)  \left(
{x}^{2}-2\,x+2 \right).
$$
Applying Theorem \ref{teojordanaditivo} we obtain, after
factorization,
$$
E(x) = -\frac{1}{2}\, \left( x-2 \right) ^{2},\quad H(x) =
-\frac{1}{2}\,{x}^{3}+ \frac{5}{2}\,{x}^{2}-4\,x+4.
$$
Applying Theorem \ref{teojordanmult} we obtain, after
factorization,
$$
e(x) = \frac{1}{4}\, \left( \sqrt {2}-2 \right)  \left(
{x}^{3}+(\sqrt {2} - 4){x}^{2}+ 4(1-\sqrt {2})x+2\,\sqrt {2} -4
\right) ,
$$
$$
h(x) =  \frac{1}{2}\, \left( \sqrt {2}-2 \right)  \left(
{x}^{3}-5\,{x}^{2}+8 \,x-8-2\,\sqrt {2} \right) ,
$$
$$
u(x) = \frac{1}{4}\,x \left( {x}^{2}-4\,x + 6\right) .
$$
One can use a mathematical software package to check that
these polynomials give the correct Jordan components of $T$.
\end{exemplo}

\section{Semisimple linear Lie algebras and groups}\label{sec2}
In this section we will obtain the Jordan decomposition in
semisimple linear Lie algebras and groups. When $\g = \gl(V)$ and
$G = \Gl(V)$, we have that $\ad(X)Y = XY - YX$  and $\Ad(g)X =
gXg^{-1}$.

\subsection{Additive Jordan decomposition}
Let $\frak{g}\subset \frak{gl}(V)$ be a semisimple Lie algebra of
$\frak{gl}(V)$.   Denote by $\n(\g)$ the normalizer of $\g$ in
$\gl(V)$, that is,
\[
\n(\g)=\{X \in \gl(V): \ad(X)\g \subset \g \}.
\]
By definition $\g$ is an ideal in $\n(\g)$.

Consider the representation $\rho: \frak{g}\rightarrow
\frak{gl}(V_{\mathbb{C}})$, of $\frak{g}$ in $V_{\mathbb{C}}$, given
by
\begin{equation}\label{eqrepresent}
\rho(X)(v)=X v,
\end{equation}
for all $X\in \frak{g}$ and $v \in V_\C$. Since $\frak{g}$ is
semisimple and $V_{\mathbb{C}}$ has finite dimension we have, by the
Weyl decomposition theorem (Theorem 3.13.1, p.222 of \cite{var}),
that there exist subspaces $V_1, \ldots, V_m$ such that
\begin{equation}\label{eqdecirred}
 V_{\mathbb{C}}=V_1\oplus \cdots \oplus V_m
\end{equation}
and each $V_k$, $1\leq k \leq m$ is invariant and irreducible by
$\rho$. For each $k=1, \ldots, m$, denote by $\g_k$ the subalgebra
of $\gl(V)$ given by
\begin{equation}\label{eqdefgk}
 \g_k=\{X \in \gl(V):\, X(V_k) \subset V_k ~{\rm and}~ \tr(X|_{V_{k}})=0\}.
\end{equation}
Since $\g$ is semisimple it follows that $\g \subset \g_k$.  Let
$\widetilde{\g}$ the following subalgebra of $\gl(V)$:
\begin{equation}\label{eqdefgtil}
 \widetilde{\g}=\n(\g)\cap \g_1 \cap \cdots \cap \g_m.
\end{equation}
We have that $\g$ is an ideal in $\wt{\g}$.

The following result is an adaptation for real semisimple Lie
algebras of proof of Theorem 6.4 , p.29 of \cite{hum}.

\begin{lema}\label{lemagtil}
With the above notations, we have that $\g = \wt{\g}$.
\end{lema}

\begin{prova}  Consider the representation $\widetilde{\rho}: \frak{g} \rightarrow
 \frak{gl}(\widetilde{\frak{g}})$ of $\frak{g}$ in
 $\widetilde{\frak{g}}$, given by $\widetilde{\rho}(X)(Y)=\ad(X)Y$, where $X \in
 \frak{g}$ and $Y \in \widetilde{\frak{g}}$. By the Weyl decomposition theorem (Theorem 3.13.1, p.222 of \cite{var}),
 there exists a subspace $\frak{h} \subset \widetilde{\frak{g}}$ invariant by $\widetilde{\rho}$
 such that
 \[
 \widetilde{\frak{g}}=\frak{g}\oplus\frak{h}.
 \]
 Since $\frak{h}$ is invariant by $\widetilde{\rho}$ and $\frak{g}$ is an ideal of
  $\widetilde{\frak{g}}$, it follows that $\ad(X)Y \in \frak{g}
 \cap \frak{h}=\{0\}$, for all $X \in \frak{g}$ and $Y
 \in \frak{h}$.

Let $Y \in \frak{h}$. For each $1\leq k \leq m$ and $v \in V_k$ we
have that
\[
\rho(X) Y(v)=XY(v)=YX(v)=Y\rho(X)v
\]
for all $X \in \frak{g}$. By the Schur lemma, there exists $c \in
\C$ such that $Y|_{V_k}=c I_k$, where $I_k$ is the identity of
$V_k$. We have that
\[
0=\tr(Y|_{V_k})= \tr(c I_k)= c \dim V_k,
\]
so that $c = 0$ and $Y|_{V_k}=0$. Since $k$ is arbitrary, $Y=0$. It
follows that $\frak{h}=\{0\}$, that is,
$\widetilde{\frak{g}}=\frak{g}$.
\end{prova}

\begin{lema}\label{lemaad}
We have the following.
 \begin{enumerate}
  \item If $E \in \frak{gl}(V)$ is elliptic then $\ad(E)$ is elliptic.
  \item If $H \in \frak{gl}(V)$ is hyperbolic then $\ad(H)$ is hyperbolic.
  \item If $N \in \gl(V)$ is nilpotent, then $\ad(N)$ is nilpotent.
 \end{enumerate}
\end{lema}

\begin{prova}
For items 1 and 2, let $\{v_1, \ldots, v_n\}$ be a basis of
$V_{\mathbb{C}}$ given by eigenvectors of a semisimple $S \in
\frak{gl}(V)$ . Let $\lambda_1, \ldots, \lambda_n$ be the respective
eigenvalues. Consider the basis of $\frak{gl}(V_{\mathbb{C}})$ given
by $E_{rs}:V_{\mathbb{C}}\rightarrow V_{\mathbb{C}}$,
$E_{rs}(v_k)=\delta_{jk}v_r$, where $\delta_{jk}$ is the Kronecker
delta.  We have that
\[
S E_{rs}v_k = S \delta_{jk}v_r =
\lambda_r\delta_{jk}v_r=\lambda_rE_{rs}v_k
\]
and thus
\[
 \ad(S)E_{rs} v_k =(\lambda_r-\lambda_s)E_{rs}v_k,
\]
which shows that $E_{rs}$ is an eigenvector of $\ad(S)$ associated
to the eigenvalue $\lambda_r-\lambda_s$.  It is then immediate that
$\ad(S)$ is elliptic (hyperbolic) when $S$ is elliptic (hyperbolic).

For the last item consider the linear map $L_N: \gl(V) \to \gl(V)$
given by $L_N(Y) = NY$, $Y \in \gl(V)$. Since $L_N^n(Y) = N^n Y$, it
follows that $L_N$ is nilpotent.  Consider also the linear map $R_N:
\gl(V) \to \gl(V)$ given by $R_N(Y) = YN$, $Y \in \gl(V)$.  Since
$R_N^n(Y) = Y N^n$, it follows that $R_N$ is also nilpotent.  Noting
that $L_N$ and $R_N$ commute, and that $\ad(N) = L_N - R_N$, it
follows that $\ad(N)$ is nilpotent.
\end{prova}

We now obtain the main result of this subsection.

\begin{teorema}
Let $\g$ be a semisimple Lie subalgebra of $\gl(V)$ and $X \in \g$.
The additive Jordan components of $X$ lie in $\g$.
\end{teorema}

\begin{prova}
Let $X = E + H + N$ be the additive Jordan decomposition of $X$.  By
Lemma \ref{lemaad} we have that $\ad(X) = \ad(E) + \ad(H) + \ad(N)$
is the additive Jordan decomposition of $\ad(X)$.  By Theorem
\ref{teojordanaditivo} it follows that $\g$ is invariant by the
Jordan components of $\ad(X)$, since they are polynomials in
$\ad(X)$.  Thus, the Jordan components of $X$ lie in $\n(\g)$. Again
by Theorem \ref{teojordanaditivo} it follows that $V_k$ is invariant
by the Jordan components of $X$, since they are polynomials in $X$.
Since $N$ is nilpotent, then $N|_{V_k}$ is also nilpotent so that
$\tr(N|_{V_k}) = 0$. Since $E$ is elliptic, then $E|_{V_k}$ is also
elliptic so that $\tr(E|_{V_k})$ is both real and pure imaginary and
thus vanishes. It follows that $\tr(H |_{V_k}) = 0$, since
$\tr(X|_{V_k}) = 0$. Hence the Jordan components of $X$ lie in
$\g_k$, for each $k=1,\ldots,m$, showing that they lie in $\wt{\g}$.
The result now follows from Lemma \ref{lemagtil}.
\end{prova}

Using the previous result and Lemma \ref{lemaad} we have the next
result, which proves also the existence of the abstract Jordan
decomposition in $\g$.

\begin{corolario}
If $\g$ is a semisimple Lie subalgebra of $\gl(V)$ then the abstract
and usual Jordan decompositions coincide.
\end{corolario}

\subsection{Multiplicative Jordan decomposition}

Let $\frak{g}$ be a semisimple Lie subalgebra of $\frak{gl}(V)$
and $G$ a connected Lie subgroup of $\Gl(V)$ with Lie algebra
$\frak{g}$.  We denote by $N(\g)$ the normalizer of $\g$ in
$\Gl(V)$ which is given by
\[
N(\g)=\{g \in \Gl(V):{\Ad}(g)\frak{g}= \frak{g}\}.
\]
We have that $G$ is a normal subgroup of $N(\g)$ since, by the
connectedness of $G$, $N(\g)$ is the normalizer of $G$ in
$\Gl(V)$.

Consider the representation $\rho: \frak{g}\rightarrow
\frak{gl}(V_{\mathbb{C}})$ given in (\ref{eqrepresent}) and the
decomposition $V_{\mathbb{C}}=V_1\oplus \cdots \oplus V_m$ given in
(\ref{eqdecirred}), such that each $V_k$ is invariant and
irreducible by $\rho$, $1\leq k \leq m$. For each $k=1, \ldots, m$,
denote by $G_k$ the subgroup of $\Gl(V)$ given by
\[
 G_k=\{g \in \Gl(V):g(V_k)\subset V_k ~{\rm e}~ \det(g|_{V_{k}})=1\}.
\]
Since $G$ is connected and semisimple it follows that $G \subset
G_k$. Consider the subgroup of $\Gl(V)$ given by
\[
 \widetilde{G}=N(\g)\cap G_1 \cap \cdots \cap G_m.
\]

\begin{lema}\label{lemaGtil}
With the above notations, $\widetilde{G}$ is algebraic and its
connected component of the identity is $G$.
\end{lema}

\begin{prova}
We first show that $\wt{G}$ is algebraic.  We start by showing that
$N(\g)$ is algebraic.  Let $\g$ be given, as a subspace, by the
kernel of $T \in \gl(\gl(V))$.  Consider $g \in \Gl(V)$, then it is
well known $\wh{g} = g^{-1} \det(g)$ is a polynomial in $g$.  Then
condition that $\Ad(g)\g = \g$ is clearly seen to be equivalent to
$T( gX\wh{g} ) = 0$, for all $X \in \g$.  Thus, taking a basis
$\{X_1, \ldots, X_n\}$ of $\g$, the condition $\Ad(g)\g = \g$ is
equivalent to the algebraic condition $T( gX_r\wh{g} ) = 0$, for
$r=1,\ldots,m$.  To show that each $G_k$ is algebraic we choose a
basis of $\{v_1,\ldots,v_n\}$ of $V_\C$ such that
$\{v_1,\ldots,v_l\}$ is a basis of $V_k$, $l \leq n$.  Denote by
$z_{rs}(g)$ the $(r,s)$-entry of the matrix of $g \in \Gl(V)$ in
this basis.  It follows that $gV_k \subset V_k$ if and only if
$z_{rs}(g) = 0$ for $r > l$, $s \leq l$, and in this case we have
$$
    \det(g|_{V_k}) = \det\left( (z_{rs}(g))_{1\leq r,s \leq l}
    \right).
$$

Now, since $G$ is connected, it is enough to show that its Lie
algebra coincides with the Lie algebra of $\widetilde{G}$.  Since
the Lie algebra of $N(\g)$ is $\n(\g)$ and the Lie algebra of $G_k$
is the subalgebra $\g_k$ given in (\ref{eqdefgk}), we have that the
Lie algebra of $\widetilde{G}$ is given by (\ref{eqdefgtil}).  The
result now follows by Lema \ref{lemagtil}.
\end{prova}

\begin{lema}\label{lemaAd}
We have the following.
 \begin{enumerate}
  \item If $e \in \Gl(V)$ is elliptic then $\Ad(e)$ is elliptic.
  \item If $h \in \Gl(V)$ is hyperbolic then $\Ad(h)$ is hyperbolic.
  \item If $u \in \Gl(V)$ is unipotent, then $\Ad(u)$ is unipotent.
 \end{enumerate}
\end{lema}

\begin{prova}
For items 1 and 2, let $\{v_1, \ldots, v_n\}$ be a basis of
$V_{\mathbb{C}}$ given by eigenvectors of a semisimple $s \in
\Gl(V)$. Let $\lambda_1, \ldots, \lambda_n$ be the respective
eigenvalues. Consider the basis of $\frak{gl}(V_{\mathbb{C}})$ given
by $E_{rs}:V_{\mathbb{C}}\rightarrow V_{\mathbb{C}}$,
$E_{rs}(v_k)=\delta_{jk}v_r$, where $\delta_{jk}$ is the Kronecker
delta.  We have that
\begin{eqnarray*}
\Ad(s) E_{rs} v_k &=& s E_{rs} s ^{-1} v_k = s E_{rs}\lambda^{-1}_k
v_k = \lambda^{-1}_k \delta_{jk} s v_r\\
&=& \lambda_r\lambda^{-1}_s\delta_{jk}v_r = \lambda_r\lambda^{-1}_s
E_{rs} v_k.
\end{eqnarray*}
which shows that $E_{rs}$ is an eigenvector of $\Ad(s)$ associated
to the eigenvalue $\lambda_r\lambda^{-1}_s$.  It is then immediate
that $\Ad(s)$ is elliptic (hyperbolic) when $s$ is elliptic
(hyperbolic).

For the last item, by Lemma IX.7.3 p.431 of \cite{helgason}, we have
that $u = {\rm e}^{N}$ where $N \in \g$ is nilpotent.  Since $\Ad(u)
= {\rm e}^{\ad(N)}$ the result follows from the last item of Lemma
\ref{lemaad}.
\end{prova}

We now obtain the principal result of this subsection.

\begin{teorema}
Let $G$ be a connected semisimple Lie subgroup of $\Gl(V)$ and $g
\in G$.  Then the multiplicative Jordan components of $g$ lie in
$G$.
\end{teorema}

\begin{prova}
Let $g = ehu$ be the multiplicative Jordan decomposition of $g$. By
Lemma \ref{lemaAd} we have that $\Ad(g) = \Ad(e)\Ad(h)\Ad(u)$ is the
multiplicative Jordan decomposition of $\Ad(g)$.  By Theorem
\ref{teojordanmult} it follows that $\g$ is invariant by the Jordan
components of $\Ad(g)$, since they are polynomials in $\Ad(g)$.
Thus, the Jordan components of $g$ lie in $N(\g)$.  Again by Theorem
\ref{teojordanmult} it follows that $V_k$ is invariant by the Jordan
components of $g$, since they are polynomials in $g$.  Since $u$ is
unipotent, then $u|_{V_k}$ is also unipotent so that $\det(u|_{V_k})
= 1$. Since $e$ is elliptic and $h$ is hyperbolic, then $e|_{V_k}$
is also elliptic and $h|_{V_k}$ is also hyperbolic so that
$\det(e|_{V_k})$ is real and have absolute value equal to one and $\det(h|_{V_k})$ is
positive real. It follows that $\det(e|_{V_k}) = \det(h|_{V_k}) =
1$, since $\det(g|_{V_k}) = 1$. Hence the Jordan components of $g$
lie in $G_k$, for each $k=1,\ldots,m$, showing that they lie in
$\wt{G}$. By Lemma \ref{lemaGtil}, it remains to show that the
Jordan components of $g$ lie in the connected component of the
identity of $\wt{G}$.

For the hyperbolic component, by Theorem \ref{teojordanmult}, we
have that $h = {\rm e}^{H}$, where $H \in \gl(V)$ is hyperbolic. Let
$\{v_1, \ldots, v_l\}$ be a basis of $V$ such that $H v_r =
\lambda_r v_r$. We have that $h^n v_r = {\rm e}^{n\lambda_r}v_r$,
for all $n \in \Z$.  In this basis, let $\{Q_s\}$ be the set of
polynomials defining $\wt{G}$.  Let $P_s$ be the polynomial obtained
by restricting $Q_s$ to the diagonal matrices in this basis.  Since
$h^n \in \wt{G}$, we have that $P_s( {\rm e}^{n\lambda_1}, \ldots,
{\rm e}^{n\lambda_l}) = 0$, for all $n \in \Z$.  By Lemma 1.142
p.116 of \cite{knapp}, we have that $P_s( {\rm e}^{t\lambda_1},
\ldots, {\rm e}^{t\lambda_l}) = 0$, for all $t \in \R$. This shows
that $e^{tH} \in \wt{G}$, for all $t \in \R$, so that $h$ lies in
the connected component of identity of $\wt{G}$.  For the unipotent
component, by Lemma IX.7.3 p.431 of \cite{helgason}, we have that $u
= {\rm e}^{N}$, where $N \in \gl(V)$ is nilpotent.  Choosing a basis
of $V$, we have that the $(r,s)$-entry of $u^n$ is a polynomial
$p_{rs}(n)$ in $n \in \Z$, since $u^n = e^{nN}$ and $N$ is
nilpotent. We have that $q_j(n) = Q_j( (p_{rs}(n))_{1 \leq r,s \leq
l})$ is also a polynomial in $n$. Since $u^n \in \wt{G}$, we have
that $q_j(n) = 0$, for all $n \in \Z$, which implies that $q_j(t) =
0$, for all $t \in \R$. This shows that $e^{tN} \in \wt{G}$, for all
$t \in \R$, so that $u$ lies in the connected component of identity
of $\wt{G}$.

Since $g$ is already in $G$ and $g = ehu$ it follows that $e$ lies
in $G$, which completes the proof.
\end{prova}

Using the previous result and Lemma \ref{lemaAd} we have the next
result, which proves also the existence of the abstract Jordan
decomposition in $G$.

\begin{corolario}
If $G$ is a connected semisimple Lie subgroup of $\Gl(V)$ then the
abstract and usual Jordan decompositions coincide.
\end{corolario}


\begin{thebibliography}{99}
\bibitem{helgason} Helgason, S.\
{\em Differential Geometry, Lie Groups and Symmetric Spaces}.
Academic Press, 1978.

\bibitem{hk} Hoffman, K.\ and Kunze, R.\ \textit{Linear Algebra}. Second Edition. Prentice-Hall, 1971.

\bibitem{hum} Humphreys, J.E.\ \textit{Introduction to Lie Algebras and Representation Theory}. Springer, 1972.

\bibitem{knapp}  Knapp, A. W.\ {\em Lie Groups Beyond an Introduction, Progress
in Mathematics}, v. 140, Birkh\"auser, 2004.

\bibitem{var} {Varadarajan, V.S.\ } \textit{Lie Groups, Lie Algebras and their Representations}.
Prentice-Hall Inc., 1974.

\bibitem{var2} {Varadarajan, V.S.\ } \textit{Harmonic Analysis on Real Reductive
Groups}. Lecture Notes in Math. \textbf{576}. Springer-Verlag, 1977.

\bibitem{warner} {Warner, G.\ } \textit{Harmonic Analysis on Semi-Simple Lie Groups I}.
Springer-Verlag, 1972.
\end{thebibliography}
\end{document}